\noindent
{\bf A.O.Prishlyak}

\noindent
{\bf CONJUGATION OF MORSE FUNCTION ON 3-MANIFOLDS} 
\medskip
The questions when two Morse function on closed manifolds are conjugated is 
investigated. Using the handle decompositions of manifolds the condition of 
conjugation is formulated. 

For each Morse function on 3-manifold the ordered generalized Heegaard diagram 
is built. The criteria of Morse function conjugation are given in the terms of 
equivalence of such diagrams.   
\bigskip
Let $M$ is a smooth manifold, $f$ and $g$ is Morse function on it. Functions 
$f$ and $g$ are called conjugated, if there exist homeomorphisms $h: M 
\rightarrow  M$, $h': R^1\rightarrow  R^1$ such, that $f\ h' = gh$ and 
homeomorphism $h'$ preserve the orientation.  

If we demand that conjugated homeomorphisms are isotopic to identity 
diffeomorphism then in the case of 1-connected $n$-manifolds with $n>5$ there 
us a criterion of Morse function conjugate in the terms of ordering basic chain
complex equivalence [1].
\medskip
\noindent
{\bf Collar handle decomposition.} 

Collar handle decomposition is a sequence of imbedding $M_0\subset M'_0 
\subset M_1\subset M'_1\subset M_2\subset ...\subset  M_N =M$ such that $M_0$ 
is a union of $n$-disks (0-handels). $M'_i$ is obtained from $M_i$ by gluing 
collar $N_i\times [0,1]$, where $N_i=\partial \ M_i$, and $M_{i+1}$ is obtained
from $M'_i$ by gluing handles. We fix the projection $\pi : N_i\times [0,1] 
\rightarrow  [0,1]$ on each collar. 

Collar handle decomposition is isomorphic, if there exist homeomorphism 
between manifolds which maps handles on handles, collars on collars and commute
with projections $\pi .$ 

Having Morse function $f: M \rightarrow  R^1$ with critical values $\{1,2,...$,
$N\}$ let us build collar handle decomposition in a such way that interiors 
of collars $N_i\times [0,1]$ will be gomeomorphic in fibers to connected 
component of $M \backslash  f^{-1}(\{1,2,...$, $N\})$.  

The first way to do this: 

Let $p$ be a critical point. Let us consider such it chart that $p=(0,0,..., 0)$ and

$$f(x) = f(p)-\sum_{i=1}^{k}{x_i^2} + \sum_{i=k+1}^{n}{x_i^2}$$ 

We denote by $D^k$, $D^{n-k}$ the disks with enough small radius 
$\hbox{$\varepsilon $}:$ 

$$D^k =\{(x_1, x_2,... , x_k,0,0,...,0):\sum_{i=1}^{k}{x_i^2} \leq \varepsilon \}$$ 

$$D^n-^k=\{(0,0,..., 0, x_{k+1},..., x_n): \sum_{i=k+1}^{n}{x_i^2} \leq \varepsilon \}.$$ 

These disks are core and cocore of handle  

$$H_j =D^k \times D^n-^k =\{(x_1, x_2,... , x_n): \sum_{i=1}^{k}{x_i^2} \leq \varepsilon ,
\sum_{i=k+1}^{n}{x_i^2} \leq \varepsilon \}.$$ 

Cut manifold $M$ by critical levels and such constructed handles for each 
critical point. We denote by $W_i$ the closure of connected components of this 
set. 

{\bf Lemma 1.} $W_i$ is gomeomorphic to $N_i\times [0,1]$, where $N_i$ is a 
regular level in $W_i.$ 

Proof. Let us fix Reiman metric on the manifold $M$ and construct vector field 
grad $f$  that is transversal to boundaries $W_i$. Because $N_i$ is a 
transversal section of grad $f$  then $W_i$ is gomeomorphic to 
$N_i\times [0,1]$.  

Consecutive gluing of handles $H_j$ and collars $N_i\times [0,1]$ we obtain 
collar handle decomposition. 

Conversely having collar handle decomposition by smoothing corners we can 
construct Morse function, which induce given handle decomposition. 

Second way: 

Cutting manifold $M$ by critical level we obtain the interiors 
$N_i\times (0,1)$ of collar $N_i\times [0,1]$. Manifold $M$ can be obtained 
under identification of corresponding points from boundaries of neighboring 
collars $N_i\times [0,1]$ and $N_{i+1}\times [0,1]$ and contraction of spheres 
$S^{k-1}$ and $S^{n-k-1}$ from $N_i\times \{1\}$ and $N_{i+1}\times \{0\}$, 
correspondingly, to critical $k$-point.  

Indeed, if we fix Reiman metric and build vector field grad $f$  then 
corresponding points from the collars boundaries is the ends of cutting 
trajectories parts. The spheres are limits sets of trajectories that correspond
to ones from stable and unstable manifolds of the critical points. 

The handle decomposition can be obtained in such way: the points from 
neighboring collars, besides points from regular neighborhood of above spheres,
are identified as previously and we glue $k$-handles to regular neighborhood of
spheres $S^{k-1}$ and $S^{n-k-1}.$ 

Reverse procedure is contraction of handles to critical points. 

{\bf Lemma 2.} Functions are conjugate if and only if correspondent collar 
handle decomposition is isomorphic.  

Proof. If there exist conjugating homeomorphism between manifolds then if we 
fix Reimann metric on one manifold we can construct induce metric on another. 
Then under this homeomorphism integral trajectories of one field will be mapped
in the integral trajectories of another. Thus from construction it is followed 
that correspondent collar handle decompositions are isomorphic. 

Let us see that collar handle decompositions don't depend from metric 
chooses. Because of all metric on manifold are homotopic then the spheres 
$S^{k-1}$ and $S^{n-k-1}$, which are constructed using its, are isotopic. Then 
these isotopies give fiber maps of collars. Extending these maps on handles we 
obtain needed homeomorphism. 

Reverse, homeomorphism between manifolds, which has collar handle 
decompositions, gives conjugated homeomorphism between Morse functions. 

It is followed from lemmas 1 and 2 that collar handle decompositions 
constructed in a first and second way are isomorphic.  
\medskip
{\bf Criteria of Morse functions conjugation.} 

Let us choose the such structure of direct product $N_i\times [0,1]$ on each 
collar that for each $t\in (0,1): N_i\times \{t\} = f ^{-1}(y)$ for
appropriate regular value $y$. Then after removing collars from manifold and 
identification $N_i\times \{0\}$ with $N_i\times \{1\}$ we obtain handle 
decomposition without collars. This decomposition depend from structure of 
direct product on each collar.  

Let us see that different structure of direct product correspond handles gluing
on isotopic embeddings, $i.e$. diffeomorphisms $p$, $q: N_i\times \{0\} 
\rightarrow  N_i\times \{1\}$ which are projection of one collar boundary 
component on another for correspondent structures of direct product are 
isotopic. Indeed, let $p_t$, $q_t: N_i\times [0,1] \rightarrow  
N_i\times \{t\}$ is projections of collar to level $N_i\times \{t\}= f 
^{-1}(y_t)$. Then isotopy can be given by formula 

$$F(x,t)= p_1(q_t(x)).$$ 

Thus the next lemma is true: 

{\bf Lemma 3.} Two Morse function are conjugate if and only if in the handle 
decomposition which associated with it correspondent handles are gluing on 
isotopic embeddings. 

Using isotopies of attaching spheres one can made handle decomposition such 
that attaching and belt spheres has transversal intersections and such that 
each handle attach to the union of less dimensional handles. We call such 
handle decomposition simple.  

The handle decomposition is called ordered if the map of handles set to set 
$\{1$, 2,..., $N\}$ is given. Each Morse function assign order in the 
corresponding handle decomposition: we can choose such $h': R^1 \rightarrow  
R^1$ that critical values set maps on the set of numbers $\{1$, 2,..., $N\}$ we
put in correspondence to each handles the number of correspondent critical 
points number.  

Ordered simple handle decompositions (OSHD) are isomorphic if there is a 
homeomorphism of the manifolds, which maps handles on handles, cores on cores, 
cocores on cocores and preserves order of handles. Denote by $M^k$ the union of
handles which indexes are no more then $k$ and $L^k = \partial \ M^k.$ 

{\bf Theorem 1.}  Two Morse functions are conjugate if and only if from first 
function OSHD we can obtain second function OSHD using 

1) attaching  spheres isotopies  in  $L^k$ with support in the boundary of handles
with less numbers; 

2) replacement handle $H_i$ by $H_i\#H_j$, if number of $H_i$ is more than the 
number of $H_j$ (handles have the same index). The new handle $H_i\#H_j$ have 
the same number as the handle $H_i.$

\medskip
Proof. Let us consider the isotopy of attaching sphere embeddings as in lemma 
3. We transform such isotopy to general position with belt spheres of less 
numbers. Then under such isotopy attaching sphere of each handle 
don't intersect belt spheres of greatest index for arbitrary 
parameter of isotopy $t$.  

If under isotopy of attaching sphere of handle $H_i$ it pass belt sphere of 
handle $H_j$ with the same index (if there is $t$ when these spheres have 
intersection) then handle $H_i$ is replaced by sum $H_i\#H_j$.  

Because such isotopy has support in the boundary of the union of handles that 
are already attached, $i.e$. handles with the less numbers, then we obtain the 
condition of theorem.   
\medskip
{\bf Generalized Heegaard diagrams} 

Let $N\cup N'=M$ be a Heegaard decomposition of the manifold $M$, 
$F=\partial \ N =\partial \ N'$ be the common surface of $N$ and $N' [2]$. The 
set $u=\{u_1$, $u_2,...$, $u_n\}$ of non intersected closed curves on surface 
$F$ are called generalized meridian system for N, if it is the boundary of 
disks $D_i \subset  H$  and if we cut  $H$ along it we obtain disconnected 
union of 3-disks. Let $v=\{v_1$, $v_2,...$, $v_m\}$ be the generalized meridian
system for $N'.$ 

The triple $(F$, $u$, $v)$ is called the generalized Heegaard diagram of 
manifold M. Diagrams $(F$, $u$, $v)$ and $(F'$, $u'$, $v' )$ are called 
gomeomorphic if there exist such homeomorphism  $h: F\rightarrow F$, that 
$h(u)=u'$, $h(v)=v'$. Diagrams $(F$, $u$, $v)$ and $(F'$, $u'$, $v' )$ are 
called semiisotopic, if there exist such isotopies $\varphi _t,\psi _t : 
F\rightarrow F$, that $\varphi _0=\psi _0 =1$, $\varphi _1(u)=u'$, 
$\psi _1(v)=v'$.  

Let us define the operation of meridian addition: the sum $u_1\# u_2$ of two 
meridian $u_1$ and $u_2$ along simple curve $\alpha $, which connect $u_1$ and 
$u_2$ is a such component of the union neighborhood boundary $\partial U(u_1 
\cup u_2\cup \alpha )$ which don't isotopic neither $u_1$ not $u_2$.

Denote by $U_1$, $U_2,...$, $U_k$ that domains which are obtained after cutting
surface $F$ along meridians $u_1$, $u_2,...$, $u_n$ and by $V_1$, $V_2,...$, 
$V_l$ {\it ---} correspondent domains for meridians $v_1$, $v_2,...$, $v_m$. 
Diagram is called ordered if the map $\sigma $ of set $\{U_1,U_2,...$, $U_k$, 
$u_1$, $u_2,...$, $u_n$, $v_1$, $v_2,...$, $v_m$, $V_1$, $V_2,...$, $V_l\}$ on 
set $\{1$, 2,..., $N\}$ is given.  

Ordered generalized Heegaard diagrams (OGHD) are called equivalent if on from 
another can be obtained using homeomorphisms, semiisotopies of diagrams (finger
moves between $u_i$ and $v_j$, if $\sigma (u_i) > \sigma (v_j))$, replacement 
meridians $u_i$ on $u_i\# u_j$ if $\sigma (u_i) < \sigma (u_j)$ and replacement
meridians $v_i$ on $v_i\# v_j$ if $\sigma (v_i) > \sigma (v_j)$. If we replace 
$u_i$ on $u_i\# u_j$ then $\sigma (u_i\# u_j) = \sigma (u_i).$  
\medskip
{\bf Morse function conjugation on 3-manifolds.} 

On 3-manifolds isotopies in 1) of theorem 1 in $L^1$ can be realized using 
finger moves or reverse moves between belt spheres of the handle $H^1$ and 
attaching sphere of the handle $H^2$ in $L^1$ under condition that number of 
the handle $H^1$ is less then number of $H^2.$ 

Indeed. Let functions are conjugate and number of the handle $H^1$ are less 
than number of $H^2$. Then when we glue handle $H^2$ belt sphere of handle 
$H^1$ and attaching sphere of handle $H^2$ can have the intersections. The 
isotopy of attaching sphere of the handle $H^2$ with respect to belt sphere of 
handle $H^1$ will be realized by finger moves and there reverse. In addition 
the pair of intersection points between these spheres (closed curves) will be 
appear or disappear.  

{\bf Theorem 2.} Two Morse function on 3-manifolds are conjugate if and only if
associated ordered generalized Heegaard diagrams are equivalent. 

Proof. Indeed, two sets of meridians from generalized Heegaard diagram are belt 
spheres of 1-handles and attaching spheres of 2-handles in $L^1$. The isotopies
from 1) of theorem 1 in $L^1$ are realized by finger moves of meridians. 1- and
2-handles additions on OGHD correspond to addition of meridians.

\bigskip
1. Sharko V.V. Functions on manifolds (algebraic and topological aspects).- Kiev: 
Naukova dumka, 1990.-196p. 

2. Matveev C.V., Fomenko A.T. Algorithmic and computer methods in 3-dimensional 
topology. --- M.: MSU, 1991.- 301p.
\vfill\eject\end